\documentclass[a4paper,10pt,leqno]{article}  
\usepackage[english]{babel}
\usepackage[centertags]{amsmath}
\usepackage{amsfonts,amssymb,amsthm}
\usepackage{graphicx}           
\usepackage[latin1]{inputenc}
\usepackage{newlfont}

\newtheorem{theorem}{Theorem}

\newtheorem{lemma}{Lemma}

\newtheorem{remark}{Remark}

\newtheorem{main result}{Main result}

\newcommand{\R}{{\mathbb R}}
\newcommand{\N}{{\mathbb N}}
\newcommand{\T}{{\mathbb T}}
\newcommand{\Z}{{\mathbb Z}}

\newcommand{\Om}{\Omega}
\newcommand{\p}{\pi}
\newcommand{\ph}{\varphi}
\newcommand{\om}{\omega}
\newcommand{\bom}{\bar{\omega}}
\newcommand{\bomm}{ ( \bom_1 , \bom_2 ) }

\newcommand{\e}{\varepsilon}
\newcommand{\m}{\mu}
\renewcommand{\a}{\alpha}
\renewcommand{\b}{\beta}
\newcommand{\s}{\sigma}

\newcommand{\g}{\gamma}
\renewcommand{\t}{\tau}
\renewcommand{\d}{\delta}

\renewcommand{\th}{\vartheta}
\newcommand{\jump}{\vspace{9pt}}   

\newcommand{\bibjump}{\vspace{-3pt}}

\newcommand{\str}{\hspace{\stretch{1}}}

\newcommand{\bs}{$\square$}

\newcommand{\lm}{\lambda}

\newcommand{\Xn}{X^{(n)}}
\newcommand{\xn}{X^{(n)}}

\newcommand{\mA}{\mathcal{A}}

\newcommand{\mR}{\mathcal{R}} 
\newcommand{\gr}{\nabla}
\newcommand{\D}{\Delta}

\newcommand{\pf}{\textbf{Proof. }}

\title{\Large{\textbf{Periodic solutions of forced Kirchhoff equations}}}

\author{\large{Pietro Baldi\footnote{Sissa, via Beirut 2--4, 34014 Trieste, Italy. \textit{E-mail:} \texttt{baldi@sissa.it}\,.}}}

\date{}

\begin{document}

\maketitle 

\vspace{-0.2cm}

\begin{abstract}
\noindent
We consider Kirchhoff equations for vibrating bodies in any dimension in presence of a time-periodic external forcing with period $2\p/\om$ and amplitude $\e$, both for Dirichlet and for space-periodic boundary conditions.

We prove existence, regularity and local uniqueness of time-periodic solutions of period $2\p/\om$ and order $\e$, by means of a Nash-Moser iteration scheme. 
The results hold for parameters $(\om,\e)$ in Cantor sets having  measure asymptotically full for $\e \to 0$.
\footnote{

\textit{2000 Mathematics Subjects Classification:} 35L70; 45K05, 35B10, 37K55.




Supported by MURST under the national project ``Variational methods and nonlinear differential equations''.}
\end{abstract}

\section{\large{Introduction}}
We consider the Kirchhoff equation 
\begin{equation}  \label{Kuni}
u_{tt} - \Delta u  \Big( 1 + \int_\Omega | \nabla u |^2 \, dx \Big) = \e g(x,t) 
 \qquad x \in \Om, \ t\in  \R 
\end{equation} 
where $g$ is a time-periodic external forcing with period $2\p/\om$ and amplitude $\e$, 
and the displacement $u:\Omega \times \R \to \R$ is the unknown. We consider both Dirichlet boundary conditions 
\begin{equation}  \label{Dirbc}
	u(x,t)=0 
\quad \forall x \in \partial\Om, \ t\in  \R
\end{equation}
where $\Omega \subset \R^d$ is a bounded, connected open set with smooth boundary,
$d \geq 1$, 
and periodic boundary conditions on $\R^d$
\begin{equation} \label{perbc}
u(x,t) = u(x+ 2\p m , t) 
\quad
\forall m \in \Z^d, \ x \in \R^d, \ t \in \R 
\end{equation}
where $\Om = (0,2\p)^d$.

Equation \eqref{Kuni} is a quasi-linear integro-partial differential equation having the structure of an infinite-dimensional Hamiltonian system, with time-depending Hamiltonian
\[
H(u,v) = \int_\Om \frac{v^2}{2} \,dx \, +  \int_\Om \frac{|\gr u|^2}{2} \,dx \, + 
\Big(  \int_\Om \frac{|\gr u|^2}{2} \,dx \Big)^2 - \int_\Om \e g u \, dx .
\]
It describes nonlinear forced vibrations of a $d$-dimensional body 
(in particular, a string for $d=1$ and a membrane for $d=2$).

This model has been proposed first in 1876 by Kirchhoff \cite{K} in dimension one, without forcing terms, with Dirichlet boundary conditions, namely
\begin{equation}  \label{1876}
	u_{tt} - u_{xx} \Big( 1 + \int_0^\p u_x^2 \, dx \Big) = 0 
	\,, \quad 	u(0,t) = u(\p,t) = 0 
\end{equation}
to describe  transversal free vibrations of a clamped string in which the dependence of the tension on the deformation cannot be neglected.
Independently, Carrier \cite{Car} and Narasimha \cite{Na} rediscovered the same equation as a nonlinear approximation of the exact model for the stretched string.

\jump

Kirchhoff equations have been studied by many authors from the point of view of the Cauchy problem
\[
u(x,0)=u_0(x) \,,
\quad
u_t(x,0)=u_1(x)
\]
starting from the pioneering paper of Bernstein \cite{Bernstein}. 
Both local and global existence has been investigated, for initial data having Sobolev or analytic regularity. 
See for example \cite{
Di-69, Pok-75, JLL-78, AS-84, DS-InvMath-92,  AP-96, Man-05} 
and the rich surveys \cite{Ar, Sp-Milano-94}.

\jump

In spite of the wide study for the Cauchy problem, to the best of our knowledge nothing is known about the existence of periodic solutions of Kirchhoff equations, except for the normal modes. 

Kirchhoff himself observed that equation \eqref{1876}, thanks to its special symmetry, possesses a sequence of normal modes, that is solutions of the form $u(x,t)=u_j(t) \sin jx$, $j=1,2,\ldots$ where $u_j(t)$ is periodic. In general, normal modes are solutions of the form $u_j(t) \ph_j(x)$ where $\ph_j(x)$ is an eigenfunction of the Laplacian on $\Om$.

In presence of a forcing term $g(x,t)$ this symmetry is broken and normal modes do not survive (except in the one-mode case $g(x,t) = g_j(t) \ph_j(x)$). 
Indeed, decomposing $u(x,t) = \sum_j u_j(t) \ph_j(x)$ shows that all components $u_j(t)$ are coupled in the integral term  $\int_\Om | \gr u |^2  dx$, and problem \eqref{Kuni} is equivalent to a system of infinitely many nonlinear coupled ODEs, namely
\[
u_j''(t) + \lm_j^2 u_j(t) \Big( 1 + \sum_k \lm_k^2 u_k^2 (t)\Big) = \e g_j(t), \ 
\quad j = 1, 2, \ldots
\]
where $g(x,t) = \sum_j g_j(t) \ph_j(x)$ and $\lm_j^2$ are the eigenvalues of the Laplacian on $\Om$.

\jump

In this paper we prove the existence of periodic solutions of \eqref{Kuni}.
We consider the amplitude $\e$ and the frequency $\om$ of the forcing term $g$ as parameters of the problem. 
We prove that there exist periodic solutions of order $\e$ and period $2\p/\om$ when $\e$ is small and $(\e,\om)$ belong to a Cantor set which has positive measure, asymptotically full for $\e \to 0$.
We prove regularity estimates for the solutions, both in Sobolev and in analytic classes, and local uniqueness (see Theorem \ref{thm:uni}, Remark  \ref{rem:both} and Theorem \ref{thm:periodic} in Section \ref{sec:results}).

\jump

There are two main difficulties in looking for periodic solutions of \eqref{Kuni}.
The first one is the so-called ``small divisors problem'', caused by resonances between the forcing frequency $\om$ with its overtones and the eigenvalues $\lm_j^2$ of the Laplacian on $\Om$.
Such  a problem arises in the inversion of the d'Alembert operator $\partial_{tt} - \D$, whose spectrum 
\[
\{ - \om^2 l^2 + \lm_j^2 \,:\, l \in \N, \ j=1,2,\ldots  \}
\]
accumulates to zero for almost every $\om$. For this reason,  $(\partial_{tt} - \D)^{-1}$ cannot map, in general, a functional space in itself, but only in a larger space of less regular functions.
This makes impossible the  application of the standard implicit function theorem.

The other difficulty is the presence of derivatives in the nonlinearity. 
In general, little is known about periodic solutions of equations of the form
\[
u_{tt} - u_{xx}  = \e f(x,t, u,u_x, u_t, u_{xx}, u_{xt}, u_{tt}) .
\]
This problem has been studied by Rabinowitz \cite{Rab} in presence of a dissipative term  $\alpha u_t$, $\a \neq 0$, and frequency $\om=1$; 
by Craig \cite{Petits} for pseudodifferential operators  
\[
u_{tt} - u_{xx} = a(x) u + b(x,|\partial_x|^{\beta} u) = 0 \,,
\quad \beta < 1 \,;
\]
by Bourgain \cite{Bou-Chicago} in cases like $u_{tt} - u_{xx} + \rho u + u_t^2 = 0$ and, for quasiperiodic solutions, \cite{Bou-IMRN-1994}  
$u_{tt} - u_{xx} = a(x) u + \e \partial_x^{1/2} (h(x,u))$.
We remark that, in general, the presence of derivatives in the nonlinearity 
makes uncertain the existence of global (even not periodic) solutions, see for example the non-existence results in \cite{Lax, Kl-Maj}  for the equation $u_{tt} - a(u_x) u_{xx} = 0$ when $a>0$, $a(v) = O(v^p)$ near $0$, $p \geq 1$. 

Our proof overcomes these two difficulties by means of a modified Newton's method in scales of Banach spaces, that is a Nash-Moser method. 
At each step of the Newton's iteration we impose some ``non-resonance conditions'' on the parameter $\om$ to control small divisors. 
For these non-resonance frequencies we can invert the linearised operator, which is a perturbation of the d'Alembertian, losing some amount of regularity. 
In this way we construct inductively a sequence of approximate solutions.  
The loss of regularity, which occurs at each step of the iteration, is 
overcome thanks to smoothing operators and to the high speed of convergence of the quadratic scheme. 

\jump

The application of Nash-Moser methods to infinite-dimensional dynamical systems having small divisors problems 
has been introduced in the Nineties by Craig, Wayne and Bourgain, in analytic or Gevray classes 
\cite{CW, Bou-Annales-98, Bou-Chicago, Petits}. 
Further developments are for example in \cite{BeBo-CantorFam, BaBe-Forced, Be-Topics, Pl-Tol}.

This technique, combined with Lyapunov-Schmidt reductions, is a flexible alternative with respect to KAM procedures \cite{W-90, Kuk-93, Poschel-Pisa}.  
In particular, currently available KAM methods seem not to apply to the quasi-linear problem \eqref{Kuni}.

\jump

Since we deal with not only analytic, but also finite order regularity, the scheme we use here differs from that in \cite{BeBo-CantorFam, BaBe-Forced} and it does not rely on analyticity assumptions. 
Such a procedure goes back directly to ideas of the original methods of \cite{Moser-61, Moser-Pisa66, Z} 
and it is developed in \cite{Be-Topics}. 
Recently \cite{BeBo-Ck} this technique has made possible to prove the existence of periodic solutions 
of nonlinear wave equations for nonlinearities having only $C^k$ differentiability.
We point out that some of the difficulties of \cite{BeBo-Ck}
are not present here, thanks to the special symmetry of the Kirchhoff nonlinearity. 
Moreover, the roles played here by space and time are inverted with respect to \cite{BeBo-Ck, BeBo-CantorFam, BaBe-Forced}.

\jump

We remark that small divisors problems become more difficult in higher dimension. For this reason, not many works deal with such problems when the dimension is larger than one (e.g. \cite{Bou-GAFA-95, Bou-Annales-98}). In that case, indeed, $\lm_j$ have a sub-linear growth, see \eqref{Weyl}. 
In general this causes further difficulties in the inversion of the linearised operators. 
In the present case, however, the structure of the Kirchhoff  nonlinear integral term makes possible the inversion in any dimension (Section \ref{sec:inv}).

\jump

Finally, we note that in case of periodic boundary conditions \eqref{perbc} zero is an eigenvalue of the Laplacian. As a consequence, we have to solve a space-average equation which is not present in the Dirichlet case \eqref{Dirbc}, see  \eqref{Y-eq} in Section \ref{sec:proof}.

\section{\large{Functional setting and main results}}  \label{sec:results}

Let $2\p / \om$ be the period of $g$. We look for solutions $u$ with the same period.
Normalising the time $t \to \om t$ and rescaling $u \to \e^{1/3} u$, \eqref{Kuni} becomes
\begin{equation} \label{K}
\om^2 u_{tt} - \Delta u = \m  \Big( \Delta u  \int_\Omega | \nabla u |^2 \, dx + g(x,t) \Big)
\end{equation}
where $\m := \e^{2/3}$ and $g,u$ are $2\p$-periodic.

\subsection{\normalsize{Case of Dirichlet boundary conditions}}
Assume that $\partial\Omega$ is $C^\infty$.
Let $ \lm_j^2, \ph_j(x) $, $j=1,2,\ldots$ be the eigenvalues and eigenfunctions of the boundary-value problem
\begin{eqnarray*}
\begin{cases}
\begin{array}{ll}
- \Delta \ph_j = \lm_j^2 \,\ph_j  & \  \text{in} \ \Om  \\
 \ph_{j }  =0 & \ \text{on} \ \partial \Om 
 \end{array}
\end{cases}
\end{eqnarray*}
with $\int_\Om \ph_j^2 \, dx = 1$ and $\lm_1 < \lm_2 \leq \ldots$\,. 
Weyl's formula for the asymptotic distribution of the eigenvalues gives  $\lm_j = O(j^{1/d})$ as $j \to \infty$, thus
\begin{equation}  \label{Weyl}
	C j^{1/d} \leq \lm_j \leq C' j^{1/d} \quad \forall j = 1,2,\ldots
\end{equation}
for some positive $C,C'$ depending on the dimension $d$ and on the domain $\Omega$
(see e.g. \cite[Vol. IV, XIII.15]{RS}).

By expansion in the basis $\{ \ph_j(x) \}$,
we define the spaces 
\begin{equation*} 
V_{\s,s}(\Om)  := \Big\{ v(x) = \sum_{j} v_j \ph_j(x) : \,  
\sum_{j} | v_j |^2 \lm_j^{2s} \, e^{2 \s \lm_j}
< \infty \Big\} 
\end{equation*}
for $s \geq 0$, $\s \geq 0$. 
Spaces $V_{0,s}$ with $\s=0$ are  used in \cite{AP-96}.
They are the domains of the fractional powers $\Delta^{s/2}$ of the Laplace operator. 
See \cite{AP-96, Fu} for a characterisation. 
For instance, $V_{0,2} = H^2(\Om) \cap H^1_0(\Om)$. 
We note that if $u \in V_{0,s}(\Om)$ then $\Delta^k u \in H^1_0(\Om) $  for all $0 \leq k \leq (s - 1)/2$. 

Spaces $V_{\s,0}$ with $s=0$ are used in \cite{AS-84}, 
where it is proved that $ \cup_{\s > 0} V_{\s,0}$ is the class of the $(-\Delta)$-analytic functions, that is, by definition,  the set of functions $v(x) \in H^1_0(\Om)$ such that
\[
\Delta^k v \in H^1_0(\Om)  
\quad \text{and} \quad
\Big| \int_\Om v \Delta^k v \, dx \Big|^{1/2} \leq C A^k k! \qquad \forall k = 0,1,\ldots
\]
for some constants $C,A$. 
In \cite{AS-84} it is observed that an important subset of $\cup_{\s > 0} V_{\s,0}$ consists of the functions $v(x)$, analytic on some neighbourhood of $\overline \Om$, such that 
\[
\Delta^k v = 0 \quad  \text{on} \ \partial \Om \qquad \forall k = 0,1,\ldots
\]
This subset coincides with the whole class of $(-\Delta)$-analytic functions when $\partial \Om$ is a real analytic manifold of dimension $(d-1)$, leaving $\Om$ on one side \cite{LM}, or when $\Om$ is a parallelepiped \cite{Ar-81}.

Clearly $V_{\s,s} = \{ u \in V_{\s,0} : \Delta^{s/2} u \in V_{\s,0} \}$ and 
$V_{\s,0} \subset V_{\s',s} \subset V_{\s',0}$ for all $s > 0$, $\s > \s' > 0$. Moreover, all finite sums $\sum_{j \leq N} v_j \ph_j(x)$ belong to $V_{\s,s}$ for all $\s,s$.

We set the problem in the spaces 
$ X_{\s,s} = H^1 ( \T, V_{\s,s})$ of $2\p$-periodic functions $u : \T \to V_{\s,s}$, 
$t \mapsto u(\cdot, t)$ with $H^1$ regularity, 
$\T := \R / 2 \p \Z$, namely
\begin{multline*} \label{def X}
X_{\s,s} := \Big\{  u(x,t) = \sum_{j \geq 1} u_j(t) \ph_j(x) : \, u_j \in H^1(\T,\R), \\
\| u \|_{\s,s}^2 := \sum_{j \geq 1} \| u_j \|_{H^1}^2 \lm_j^{2s} \, e^{2 \s \lm_j}
< \infty \Big\} .
\end{multline*}

\begin{theorem}  \label{thm:uni} 
\emph{(Case of Dirichlet boundary conditions).}
Suppose that $g \in X_{\s,s_0}$ for some $\s \geq 0$, $s_0 > 2d$. 
Let $s_1 \in (1 +d, 1 + s_0/2)$.
There exist positive constants $\d,C$ with the following properties.

For every $\g \in (0,\lm_1)$ there exists a Cantor set $\mA_\g \subset (0,+\infty) \times (0,\d \g)$ of parameters such that for every $(\om,\m) \in \mA_\g$ there exists a classical solution $u(\om,\m) \in X_{\s, s_1}$ of \eqref{K}\eqref{Dirbc}.
Such a solution satisfies 
\begin{equation*}  
\| u(\om,\m) \|_{\s, s_1} \leq \frac{\m}{\g}\, C ,
\qquad
\| u(\om,\m)_{tt} \|_{\s, s_1 -2} \leq \frac{\m}{\g \om^2} \, C 
\end{equation*}
and it is unique in the ball $\{ \| u \|_{\s, s_1} < 1 \}$.

The set $\mA_\g $ satisfies the following Lebesgue measure property: for every $0<\bom_1 < \bom_2 < \infty$ there exists a constant $\bar C$ independent on $\g$  
such that in the rectangular region $\mR_\g := \bomm \times (0,\d \g)$ there holds
\[
\frac{ | \mR_\g \cap \mA_\g | }{ | \mR_\g | } \, > 1 - \bar C \, \g \,. 
\]
\end{theorem}

\smallskip

We recall that \eqref{K} is obtained from \eqref{Kuni} 
by the normalisation $t \to \om t$ and the rescaling $u \to \e^{1/3}u$.
Hence, going back, the solution $u(\om,\mu)$ of \eqref{K} found in Theorem \ref{thm:uni} 
gives a solution of \eqref{Kuni} of order $\e$ and period $2\p/\om$.

\smallskip

\begin{remark}  \label{rem:both}
\emph{Theorem \ref{thm:uni} covers both Sobolev and analytic cases:}
\begin{itemize}
\item
(Sobolev regularity). 
\emph{If $g$ belongs to the Sobolev space $X_{0,s_0}$,
then the solution $u$ found in the theorem belongs to the Sobolev space $X_{0,s_1}$.} 

\item
(Analytic regularity). \emph{If $g$ belongs to the analytic space $X_{\s_0,0}$, then 
$g \in X_{\s_1, s_0}$ for all $\s_1 \in (0,\s_0)$. Indeed,}
\[
\frac{ \xi^{s_0}}{\exp[(\s_0 - \s_1) \xi]} \,
\leq \Big( \frac{s_0}{(\s_0 - \s_1) e} \, \Big)^{s_0} =: C
\quad \forall \xi \geq 0 \,,
\]
\emph{therefore} 
\[
\| g \|_{\s_1, s_0}^2 = \sum_j \| g_j \|_{H^1}^2 \lm_j^{2 s_0} e^{2 \s_1 \lm_j} \, \frac{e^{2 \s_0 \lm_j}}{e^{2 \s_0 \lm_j}} \,
\leq C^2\| g \|_{\s_0,0}^2 \,.
\]
\emph{Since $g \in X_{\s_1, s_0}$, the solution $u$ found in the theorem belongs to the analytic space $X_{\s_1,s_1} \subset X_{\s_1,0}$.}
\end{itemize}
\end{remark}

\begin{remark} \label{rem:bootstrap}
\emph{If $g (x,\cdot) \in H^{r}(\T)$, $r \geq 1$, then the solution $u$ of \eqref{Kuni} found in the theorem satisfies $u(x,\cdot) \in H^{r+2}(\T)$ by bootstrap.} 
\end{remark}

\begin{remark}  \label{rem:3D}
(Nonplanar vibrations).
\emph{We can consider the Kirchhoff equation for a string in the 3-dimensional space 
\begin{equation} \label{3D}
u_{tt} - u_{xx} \Big( 1 + \int_0^\p |u_x|^2 \, dx \Big) = \e g(x,t),
\qquad
g = \Big( \begin{array}{c} \! g_1 \! \\ \! g_2 \! \end{array}\Big),
\quad 
u = \Big( \begin{array}{c} \! u_1 \! \\ \! u_2 \! \end{array} \Big)
\end{equation}
where the forcing $g$ and the displacement $u $ are $\R^2$-vectors belonging to the plane orthogonal to the rest position of the string, see \cite{Car,Na}. 
In this case nonplanar vibrations of the string are permitted.}

\emph{Setting $\| u_j \|_{H^1}^2 := \| u_{1,j} \|_{H^1}^2 + \| u_{2,j} \|_{H^1}^2$ in the definition of the spaces $X_{\s,s}$, Theorem \ref{thm:uni} holds true for problem \eqref{3D} as well.}
\end{remark}

\subsection{\normalsize{Case of periodic boundary conditions}}
The eigenvalues and eigenfunctions of the Laplacian on $\T^d$ are $|m|^2$, $e^{im \cdot x}$ for $m \in \Z^d$.
We consider a bijective numbering $\{ m_j : j \in \N \}$ of  $ \Z^d $ 
such that $|m_j| \leq |m_{j+1}|$ for all $j \in \N = \{ 0, 1, \ldots \}$, and we denote
\[
\tilde \lm_j^2 := | m_j |^2 , 
\quad
\tilde \ph_j(x) := e^{im_j \cdot x}
\quad \ 
\forall j \in \N .
\]
We note that $\tilde \lm_0 =0$, $\tilde\ph_0(x) \equiv 1 $ and $\tilde \lm_j \geq 1$ for all $j \geq 1$.
Weyl's estimate \eqref{Weyl} holds true for $\tilde\lm_j$ as well, because the number of integer vectors $m \in \Z^d$ such that $|m| \leq \lm$ is $O(\lm^d)$ for $\lm \to +\infty $ (see \cite[Vol. IV, XIII.15]{RS}).

We define
\begin{multline*} \label{def X tilde}
\tilde X_{\s,s} := \Big\{  u(x,t) = \sum_{j \geq 0} u_j(t) \tilde \ph_j(x) : \, u_j \in H^1(\T,\R), \\
\| u \|_{\s,s}^2 := \| u_0 \|_{H^1}^2 + \sum_{j \geq 1} \| u_j \|_{H^1}^2 \tilde \lm_j^{2s} \, e^{2 \s \tilde \lm_j} < \infty \Big\} .
\end{multline*}

\begin{theorem}  \label{thm:periodic}
\emph{(Case of periodic boundary conditions).}
Suppose that $g \in \tilde X_{\s,s_0}$ for some $\s \geq 0$, $s_0 > 2d$, and 
\begin{equation}  \label{hyp:ort}
\int_{(0,2\p)^{d+1}} g(x,t) \, dx dt = 0 \,.
\end{equation}
Let $s_1 \in (1 +d, 1 + s_0/2)$.
There exist positive constants $\d,C$ with the following properties.

For every $\g \in (0,1)$ there exists a Cantor set $\mA_\g \subset (0,+\infty) \times (0,\d \g)$ of parameters such that for every $(\om,\m) \in \mA_\g$ there exists a classical solution $u(\om,\m) \in \tilde X_{\s, s_1}$ of \eqref{K}\eqref{perbc} satisfying 
\[
\int_{(0,2\p)^{d+1}} u(\om,\m)(x,t) \, dx dt = 0 \,.
\]
Such a solution satisfies  
\begin{equation} \label{estim sol per} 
	\| u(\om,\m) \|_{\s,s_1} \leq \frac\m\g \, \Big( 1 + \frac{1}{\om^2 } \Big) \, C \,,
\qquad
\| u(\om,\m)_{tt} \|_{\s,s_1-2} \leq \frac{\m}{\g \om^2} \, C 
\end{equation}
and it is unique in the ball 
$\{ \int_{(0,2\p)^{d+1}} u (x,t) \, dx dt = 0 , \ \| u \|_{\s, s_1} < 1 \}$. 

\vspace{1pt}

The set $\mA_\g $ satisfies the same measure property of Theorem \ref{thm:uni}.
\end{theorem}

\begin{remark} \label{rem:modulo costanti}
\emph{If $u(\om,\m)$ is a solution of \eqref{K}\eqref{perbc}, then also $u(\om,\m) + c$, $c \in \R$, solves \eqref{K}\eqref{perbc}.}
\end{remark}

\subsection{\normalsize{Outline of the proof}}

The rest of the paper is devoted to the proof of the theorems. 
In Sections \ref{sec:iteration},\ref{sec:solution},\ref{sec:Cantor},\ref{sec:inv} we develop the details for the proof of Theorem \ref{thm:uni}, then the same calculations are used to prove Theorem \ref{thm:periodic} in Section \ref{sec:proof}.

In Section \ref{sec:iteration} we perform the Nash-Moser iteration to construct the approximating sequence $(u_n)$, for $\m$ small and $(\om,\m)$ belonging to smaller and smaller ``non-resonant'' sets $A_n$. Avoiding resonances allows to invert the linearised operator at each step of the iteration.

In Section \ref{sec:solution} we prove that $u_n$ converges to a solution of the Kirchhoff equation if $(\om,\m) \in A_n$ for all $n$. Local uniqueness of the solution is also proved.

In Section \ref{sec:Cantor} we prove that the intersection of all $A_n$ is a nonempty set, which is very large in a Lebesgue measure sense.

In Section \ref{sec:inv} we prove the invertibility of the linearised operator for $(\om,\m) \in A_n$ and we give an estimate on the inverse operator.

In Section \ref{sec:proof} we complete the proof of Theorem \ref{thm:uni} and we prove Theorem \ref{thm:periodic}.

\section{\large{The iteration scheme}}  \label{sec:iteration}
We fix $\s \geq 0$ once for all. 
In the following, we write in short $X_s := X_{\s,s}$, $\| u \|_s := \| u \|_{\s,s}$. 
We remark that all the following calculations  holds true both in the Sobolev case $\s=0$ and in the analytic case $\s >0$. 
Indeed, the only index used in the present Nash-Moser method is $s$.

We set the iterative scheme in the Banach spaces $X_s$ endowed with the smoothing operators $P_n$, defined in the following way. 
We consider a constant $ \chi \in (1,2)$ and denote 
\begin{equation} \label{Nn}
N_n := \exp( \chi^n) 
\end{equation}
for all $n \in \N$.  
We define the finite-dimensional space 
\[
\xn := \Big\{ u(x,t)= \sum_{\lm_j \leq  N_n} u_j(t) \ph_j(x) \Big\}
\]
and indicate $P_n$ the projector onto $ \xn$ (truncation operator).
For all $s,\a \geq 0$ there holds the smoothing properties
\begin{eqnarray} 
	 \| P_n u \|_{s+\a} & \leq & N_n^\a \| u \|_s
	\quad \forall u \in X_s  
	\label{S1} \\
	 \| (I-P_n)u \|_s  & \leq & N_n^{-\a} \, \| u \|_{s+\a} 
	\quad \forall u \in X_{s+\a} 
	\label{S2} 
\end{eqnarray}
where $I$ is the identity map.
We denote
\begin{eqnarray*}
&L_\om := \om^2 \partial_{tt} - \Delta \,,
\quad
f(u) := \Delta u \int_\Omega | \nabla u|^2 \, dx, & \\
& F(u) := L_\om u - \m f(u) - \m g &
\end{eqnarray*}
so that \eqref{K} can be written as
\begin{equation} \label{Kf}
F(u)=0 \,.
\end{equation}
Note that $f$ is not a composition operator, because of the presence of the integral. 
The map $f$ is cubic: indeed $f(u)=A[u,u,u]$ where $A$ is the three-linear map 
$A[u,v,w]= \Delta u \int_\Om \nabla v \circ \nabla w$ $dx$. 
Moreover, since the integral term $\int_\Om | \nabla u|^2 dx$ depends only on time, 
there holds
\begin{equation*} 
	f(u) \in \xn \quad \forall u \in \Xn .
\end{equation*}
The quadratic remainder of $f$ at $u$ is
\begin{eqnarray} \label{Q}
Q (u,h) 
& := &  f(u+h) - f(u) - f'(u)[h] \\
& = & \D u \int_\Om | \gr h |^2 dx + 
\D h \, \int_\Om (2 \gr u \circ \gr h + | \gr h |^2 )\, dx . \nonumber
\end{eqnarray}
We observe that, if $a(t)$ depends only on time, then 
\[
\| a(t) u(x,t) \|_s \leq \| a \|_{H^1} \| u \|_s 
\]
(we omit a factor given by the algebra constant of $H^1(\T)$). As a consequence, by H\"older inequality it is easy to estimate $\| f(u) \|_s, \| f'(u)[h]\|_s$ and $\| Q(u,h) \|_s$.

\bigskip

We adapt the Newton's scheme with smoothing operators $P_n$ to the special structure of problem \eqref{K}.  
We will construct a sequence $(u_n)$ defining 
\begin{eqnarray} \label{scheme}
	u_0 := 0, 
	\quad
	u_{n+1} :=  u_n - F'(u_n)^{-1} [ L_\om u_n - \m f(u_n) - \m P_{n+1} g] 
\end{eqnarray}
provided the linearised operator 
\[
F'(u_n) : h \mapsto L_\om h - \mu f'(u_n)[h]
\]
admits a bounded inverse $F'(u_n)^{-1}$ on $X^{(n+1)}$. 
In this inversion problem a small divisors difficulty arises. 
We will prove (Lemma \ref{lemma:inv}) that $F'(u_n)$ can be inverted if the parameters $(\om,\m)$ belong to some ``nonresonant'' set $A_{n+1}$, defined as follows.
First, 
\[
A_0 := (0, +\infty) \times (0,1) .
\]
By induction, suppose we have constructed $A_{n}$ and $u_n$. 
We denote
\begin{equation}  \label{an}
	a_n(t) := \int_\Om | \gr u_n|^2 \, dx ,
\end{equation}
we consider the Hill's eigenvalue problem 
\begin{equation*}  
\begin{cases}
 y'' + p^2 \big( 1 + \m a_n(t) \big) \,y = 0  \\
 y(t) = y(t + 2\p) 
\end{cases}	
\end{equation*}
and indicate $(p_{l}^{(n)})^2$ its eigenvalues, $l \in \N $.
For 
\[
\t >d, \quad \g \in (0,\lm_1)
\]
we define
\begin{equation}  \label{An+1}
A_{n+1} := \Big\{ (\om,\m) \in A_{n} \, : \, 
| \om p_l^{(n)} - \lm_j| > \frac{\g }{\lm_j^{\t}} \quad \forall \lm_j \leq N_{n+1}, 
\ \  l \in \N \Big\} \,.
\end{equation}
\begin{remark} \label{rem:open}
\emph{Note that for all $\m,n$ the set $A_n(\m) := \{ \om : (\om,\m) \in A_n \}$ is open. Indeed, for every $0 < \bar \om_1 < \bom_2 < \infty$ the intersection $(\bom_1, \bom_2) \cap A_n(\mu)$ is defined by means of finitely many strict inequalities (see \eqref{Weyl} and \eqref{comparison lmm}).}
\end{remark}

We fix a positive constant $R$ such that, 
if $u \in X_{1}$ and $\| u \|_{1} < R$, 
then $a(t) := \int_\Om | \gr u|^2 \, dx$ satisfies 
$\| a \|_{H^1} < 1$ and $\| a \|_{\infty} < 1/2$. 

\begin{lemma}  \label{lemma:inv} 
\textnormal{(Inversion of the linearised operator). } 
There exist two universal constants $K_1,K_1'$ with the following property.
Let $u  \in X^{(n)}$ with  $\| u \|_{1} < R$.
Let $(\om,\m) \in A_{n+1}$. If
\begin{equation}  \label{eccesso}
 \frac{\m}{\g  } \, \| u \|_{\t+1}^2  < K_1' \, ,	
\end{equation}
then $F'(u )$ is invertible, $F'(u )^{-1} : X^{(n+1)} \to X^{(n+1)}$ and
\begin{equation}  \label{stima inverso}
\| 	F'(u )^{-1} h \|_{0} \leq 
\frac{K_1 }{\g } \, 
\| h \|_{\t-1} 
\quad
\forall h \in X^{(n+1)} \,.
\end{equation}
\end{lemma}

\noindent
\pf  
In Section \ref{sec:inv}. 
\str \bs

\begin{lemma}  \label{lemma:indutt}
\textnormal{(Construction of the approximating sequence).}
Let $g \in Y_{s_0}$, $s_0 > 2d$.    
Let $\t \in ( d, s_0/2)$. 
There exist a choice for $\chi$ in the definition \eqref{Nn} and  positive constants $K,b, \d_0$, with 
$b (2 - \chi) > \t+1$, satisfying the following properties.

\smallskip

\noindent
\emph{(First step).} If $(\om,\m) \in A_1$ and $\m/\g < \d_0$, then there exists $u_1 \in X^{(0)}$ defined by \eqref{scheme}, and there holds 
\begin{equation}  \label{u1}
\| u_1 \|_0 < K \frac{\mu}{\g} \, \exp(- b \chi) \,.
\end{equation}

\noindent
\emph{(Induction step).} Suppose we have constructed $u_1, \ldots, u_n$
by \eqref{scheme} for $(\om,\m) \in A_{n}$, $n \geq 1$,
where each $A_{k+1}$ is defined  by means of $u_{k}$ by \eqref{An+1}, 
and $u_k \in X^{(k)}$. Suppose that $\m/\g < \d_0$. Let 
\[
h_{k+1} := u_{k+1} - u_k.
\] 
Suppose that for all $k=1, \ldots, n$ there holds 
\begin{equation}  \label{hk}
	\| h_k \|_0 < K \frac{\mu}{\g} \, \exp(- b \chi^k) \,.
\end{equation}
If $(\om,\m) \in A_{n+1}$ then there exists $h_{n+1} \in X^{(n+1)}$ defined by \eqref{scheme} and there holds 
\begin{equation} \label{hn+1} 
	\| h_{n+1} \|_0 < K \frac{\mu}{\g} \, \exp(- b \chi^{n+1}) \,.
\end{equation}
\end{lemma}

\noindent
\pf
\emph{(First step).} Since $u_0=0$ and $(\om,\m) \in A_1$, by Lemma \ref{lemma:inv} $F'(0)$ is invertible and \eqref{scheme} defines 
\begin{equation}  \label{u1 def}
u_1 = - F'(0)^{-1} [F(0) + \m (I - P_1)g] 
= \m L_\om^{-1} P_1 g \,.
\end{equation}
By \eqref{stima inverso}, the inequality \eqref{u1} holds true provided
\begin{equation} \label{prov 2}
K_1 \| g \|_{s_0} < K \exp(-b\chi), \quad \t-1 \leq s_0 \,.
\end{equation}

\smallskip

\emph{(Induction step).}
To define $h_{n+1}$ by \eqref{scheme}, we have to verify the hypotheses of Lemma \ref{lemma:inv}. 
By \eqref{S1} and \eqref{hk} 
\begin{equation}  \label{hk t+1}
 \| h_k \|_{\t+1} \leq N_k^{\t+1} \| h_k \|_0 
\leq K \, \frac{ \m}{\g} \, \exp[ (- b+\t+1) \chi^k ] 
\end{equation}
because $h_k \in X^{(k)}$. Then 
\begin{equation} \label{sommina}
	\| u_n \|_{\t+1} \leq \sum_{k=1}^n \| h_k \|_{\t+1} 
	< K \frac{\mu}{\g} \, \sum_{k=1}^{+\infty} \exp[(- b +\t+1) \chi^k] 
\end{equation}
which is finite for $b > \t+1$. Thus condition \eqref{eccesso} is verified provided
\begin{equation}  \label{prov 3}
b > \t+1,
\quad
K^2 \, \Big( \frac{\m}{\g} \Big)^3 \, 
C_0^2  < K_1',
\quad
C_0 := \sum_{k \geq 1} \exp[(- b +\t+1) \chi^k]\,.
\end{equation}
Since $\| u_n \|_1 \leq \| u_n \|_{\t+1}$, by \eqref{sommina} we have $\| u_n \|_1 <R$ provided
\begin{equation}  \label{prov 3.5}
K \, \frac{\m}{\g} \, C_0 < R \,.
\end{equation} 

Since $(\om,\m) \in A_{n+1}$, we can apply Lemma \ref{lemma:inv} and we define $h_{n+1}$ according to the scheme, namely
\begin{equation}  \label{precomodo}
h_{n+1} := 
- F'(u_n)^{-1} [ F(u_n) + \m (I-P_{n+1})g] \,.
\end{equation}
By \eqref{stima inverso} we have
\begin{equation}  \label{intermedio}
\| h_{n+1} \|_0 \leq \frac{K_1 }{\g } \, \| F(u_n) + \m (I-P_{n+1})g \|_{\t-1} \,.
\end{equation}
By construction \eqref{scheme}, $u_{n}$ satisfies 
\begin{equation*}  
F'(u_{n-1}) h_n = - L_\om u_{n-1} + \m f(u_{n-1}) + \m P_n g .
\end{equation*}
By Taylor expansion $F(u_n) = F(u_{n-1}) + F'(u_{n-1}) h_n - \mu Q(u_{n-1}, h_n)$, where $Q$ is defined in \eqref{Q}. Thus
\begin{eqnarray} \label{comodo}
F(u_n) =  -\mu [ (I-P_n) g +  Q(u_{n-1}, h_n)]
\end{eqnarray}
and \eqref{intermedio} gives
\begin{equation}  \label{hPQ}
\| h_{n+1} \|_0 \leq \frac{K_1 \m }{\g } \,
\| (P_{n+1} - P_n ) g + Q(u_{n-1}, h_n) \|_{ \t-1} \,.
\end{equation}
Now $(P_{n+1} - P_n)g = (I-P_n)P_{n+1}g$, then by \eqref{S2} 
\[
\| (P_{n+1} - P_n )   g \|_{\t-1} \leq 
\frac{1}{N_{n}^\b} \, \| P_{n+1} g \|_{ \t-1 + \b} 
\leq \frac{1}{N_{n}^\b} \, \| g \|_{s_0} 
\]
for 
\begin{equation*}  
\t-1 + \b \leq s_0 , \quad 
\b >0 \,.
\end{equation*}
To estimate $\| Q(	u_{n-1}, h_n) \|_{ \t-1}$, we note that 
\[
\| \D u_{n-1} \int_\Om | \gr h_n |^2  \|_{\t-1} 
\, \leq \, \| u_{n-1} \|_{\t+1} \| h_n \|_1^2 
\, < \, K \, \frac\m\g \, C_0 \, \| h_n \|_1^2
\]
by \eqref{hk t+1}, and $\| h_n \|_1 \leq N_n \| h_n \|_0$ by \eqref{S1}.
For the second term, recalling that $2 < \t+1$,
\begin{eqnarray*}
\| \D h_n \int_\Om \gr (2 u_{n-1} + h_n ) \circ \gr h_n \|_{\t-1} 
& \leq & \| h_n \|_{\t+1} \| h_n \|_0 \| 2 u_{n-1} + h_n \|_2  \\
& < & 2 K \, \frac\m\g \, C_0 \, \| h_n \|_{\t+1} \| h_n \|_0
\end{eqnarray*}
by \eqref{hk t+1}, 
and $\| h_n \|_{\t+1} \leq N_n^{\t+1} \| h_n \|_0$ by \eqref{S1}.
Then 
\[
\| Q(u_{n-1}, h_n) \|_{ \t-1} 
< 3 K \, \frac \m\g \, C_0 N_n^{\t+1} \| h_n \|_0^2 \,.
\]
As a consequence,  \eqref{hn+1} holds true provided
\begin{equation} \label{prov 4}
	K_1 \, \frac{1}{N_{n}^\b} \, \| g \|_{s_0} 
	< \frac12 \, K \exp(-b\chi^{n+1})
\end{equation}
and
\begin{equation} \label{prov 5}
		3 K_1 \, \frac\m\g \, C_0 \, N_n^{\t+1} \| h_n \|_0^2 
		< \frac12 \,  \exp(-b\chi^{n+1}) \,.
\end{equation}
Condition \eqref{prov 4} is satisfied for 
\begin{equation} \label{prov 6}
	 \b > b \chi, \qquad 
	K > \frac{2 K_1 \| g \|_{s_0} }{\exp[( \b - b \chi) \chi]} 
\end{equation}
and, by \eqref{hk}, condition \eqref{prov 5} is satisfied for 
\begin{equation}  \label{prov 7}
	b (2 - \chi) > \t+1 , \qquad
	\frac{\m}{\g} < 
\Big\{ \frac{\exp[ ( b(2-\chi) -\t-1) \chi]}{6 K_1 C_0 K^2} \Big\}^{1/3} .
\end{equation}
Since $2\t < s_0$, we can fix $\chi \in (1,2)$ so close to 1 that 
\begin{equation*}  
 \t - 1 + (\t +1) \frac{\chi}{2 - \chi} \, < s_0\,.
\end{equation*}
Now we fix $b$ such that 
\[
b (2 - \chi) > \t+1  \, , 
\qquad 
 \t -1 + b \chi < s_0 
\]
and then we fix $\b$ as
\[
\b = s_0 -\t +1 \,.
\]
So \eqref{prov 6} and \eqref{prov 2} are satisfied for $K$ big enough, and we fix $K$ in such a way.
Then \eqref{prov 3},\eqref{prov 3.5} and \eqref{prov 7} are satisfied for $\m/\g$ small enough. 
\str \bs

\section{\large{The solution}}  \label{sec:solution}

\begin{lemma}  \textnormal{(Existence of a solution).}  \label{lemma:existence}
Assume the hypotheses of Lemma \ref{lemma:indutt} and suppose that $(\om,\e) \in A_n$ for all $n \in \N$.  
Then the sequence $(u_n)$ constructed in Lemma \ref{lemma:indutt} converges in 
$X_{\t +1}$ to $ u_\infty 	:= \sum_{k \geq 1} h_k$.\, 
$u_\infty$ is a solution of \eqref{Kf} and 
\begin{equation}  \label{stima u infty}
\| u_\infty \|_{\t+1} \leq \frac{\m}{\g } \, C
\end{equation}
for some $C$. 
Moreover, $(u_n)_{tt}$ converges to $(u_\infty)_{tt}$ in $X_{\t-1}$, 
\begin{equation}  \label{H3}
\| (u_\infty)_{tt} \|_{\t-1} \leq \frac{\m}{\g \om^2} \, C
\end{equation}
so $ u_\infty $ is a classical solution of \eqref{K}.
\end{lemma}

\noindent
\pf 
By \eqref{hk t+1}, the series $\sum_k \| h_k \|_{\t+1}$ converges, $u_n$ converges to $u_\infty$ in $X_{\t+ 1}$ 
and \eqref{stima u infty} holds true.

By \eqref{u1 def} there holds  $\om^2 (h_1)_{tt} = \D h_1 + \m P_0 g$.  
 By \eqref{S1} and \eqref{hk}
\[
\| \D h_1 \|_{\t-1} + \m \| P_0 g \|_{\t-1} 
\leq \frac{K \m}{\g} \, \exp[ (- b+\t+1) \chi] + \m \| g \|_{s_0} ,
\]
so that 
\begin{equation}  \label{h1 tt}
 \| (h_1)_{tt} \|_{\t-1} \leq C \, \frac{\m}{\g \om^2 }
\end{equation}
for some $C$ (recall that $\g < \lm_1$).
For $n \geq 1$, by \eqref{precomodo} and \eqref{comodo}
\begin{equation}  \label{comodissimo}
F'(u_n) h_{n+1} = \m [ (P_{n+1} - P_{n})g + Q(u_{n-1},h_n)] ,
\end{equation}
thus
\[
\om^2 (h_{n+1})_{tt} = \D h_{n+1} + \m 
\big( f'(u_n)[h_{n+1}] + (P_n - P_{n-1})g + Q(u_{n-1},h_n) \big) .
\]
By \eqref{hk t+1},\eqref{hPQ},\eqref{prov 4} and \eqref{prov 5} we get
\begin{equation}  \label{hk tt}
 \| (h_{n+1})_{tt} \|_{ \t-1} \leq C \,\frac{\m}{\g \om^2} \, 
 \exp[ (- b+\t+1) \chi^{n+1}] .
\end{equation}
It follows that $(u_n)_{tt}$ converges in $X_{\t-1}$,\, $(u_\infty)_{tt} \in X_{\t-1}$, so that $u_\infty$ has regularity $H^3 \subset C^2$ in time and \eqref{H3} holds true. 
As a consequence $F(u_n)$ converges to $ F(u_\infty)$ in $X_{\t-1}$.

On the other hand, by \eqref{comodo},\eqref{prov 4} and \eqref{prov 5}
\begin{equation*}  
\| F(u_n) \|_{\t-1} < \frac{K \mu}{K_1} \,   \exp(-b \chi^{n+1}),
\end{equation*}
and $F(u_n) \to 0$ in $X_{\t-1}$ . 
Then $F(u_\infty)=0$.
\str \bs

\smallskip

\begin{remark}  \label{rem:niente paura}
\emph{We will prove in Lemma \ref{lemma:Cantor} that the set $\{ (\om,\m) \in A_n \ \forall n \in \N \}$ is nonempty and has positive, large measure. 
As a consequence, the sequence $(u_n)$ of Lemma \ref{lemma:existence} is defined for all $(\om,\m)$ in that large set.}
\end{remark}

\begin{lemma}  
\textnormal{(Uniqueness of the solution).}  \label{lemma:uniqueness}
Assume the hypotheses of Lemma \ref{lemma:existence}. 
There exists $\d_1 \in (0,\d_0]$ such that,
for $\m / \g < \d_1 $,  $u_\infty$ is the unique solution of \eqref{Kf} in the ball
$\{ v \in X_{\t+1} : \| v \|_{\t+1} < 1 \}$.
\end{lemma}

\noindent
\pf
Suppose $v $  is another solution of \eqref{Kf}, 
with $\| v \|_{\t+1} < 1$.  
Let $v_n := P_n v$. 
Projecting the equation $F(v) = 0 $ on $X^{(n)}$ gives
\begin{equation*}
	L_\om v_n = \m ( f(v_n) + R_n(v) + P_n g ) ,
	\quad \,
	R_n(v) := \D v_n \int_\Om | \gr(v - v_n)|^2 \, dx  .
\end{equation*}
Since $u_n$ solves \eqref{scheme}, that is
\[
F'(u_n) h_{n+1} = - L_\om u_n + \m (f(u_n) + P_{n+1} g ),
\]
the difference  $w_n := v_n - u_n$ satisfies
\[
 L_\om w_n - \m (f(v_n) - f(u_n)) - \m R_n(v) + \m (P_{n+1} - P_n) g 
 = F'(u_n) h_{n+1} \,.
\]
Since $f(v_n) - f(u_n) = f'(u_n)[w_n] + Q(u_n , w_n)$, applying $F'(u_n)^{-1}$ 
\begin{equation}  \label{frutta}
w_n = h_{n+1} + \m F'(u_n)^{-1} [ Q(u_n,w_n) + R_n(v) - (P_{n+1} - P_n) g] \,.
\end{equation}
Now, by \eqref{stima inverso} 
\[
\| \m F'(u_n)^{-1} Q(u_n,w_n) \|_{0} \leq K_1 \, \frac\m\g \, \| Q(u_n,w_n) \|_{\t-1} \,.
\]
By assumption and \eqref{stima u infty}
\[
\| w_n \|_2 \leq \| w_n \|_{\t+1} \leq \| v \|_{\t+1} + \| u_\infty \|_{\t+1} < C
\]
and by \eqref{hk t+1} $ \| u_n \|_{2} \leq \| u_n \|_{\t+1} < C$, so that
$\| Q(u_n,w_n) \|_{\t-1} \leq C \| w_n \|_0$ and
\begin{eqnarray*}
\| \m F'(u_n)^{-1} Q(u_n,w_n) \|_{0} \leq K_1 \, \frac\m\g \, C \, \| w_n \|_0
\leq \frac12 \, \| w_n \|_0
\end{eqnarray*}
provided $\m / \g$ is small enough. Thus \eqref{frutta} gives
\[
\frac12 \, \| w_n \|_0 \leq \| h_{n+1} \|_0
+ K_1 \, \frac\m\g \, \| R_n(v) - (P_{n+1} - P_n) g \|_{\t-1} \,.
\]
By \eqref{hn+1} and \eqref{S2} the right-hand side tends to 0 as $n \to \infty$, 
so that $\| v_n - u_n \|_{0} \to 0$. 
Since $v_n$ converges to $v$ and $u_n $ to $u_\infty$ in $X_0$, 
it follows that $v = u_\infty$.
\str \bs

\section{\large{The Cantor set of parameters}}  \label{sec:Cantor}

\begin{lemma}  \label{lemma:der un}
\textnormal{(Regular dependence on the parameter $\om$).} 
Assume the hypotheses of Lemma \ref{lemma:indutt}.
There exist $\d_2 \in (0, \d_1]$ such that all the maps 
\[
h_n \, : \, A_n \cap \{ (\om,\m) : \m / \g < \d_2 \} \to X^{(n)} \, , 
\quad 
(\om,\m) \mapsto h_n(\om,\m) 
\]
are differentiable w.r.t. $\om$ and
\begin{equation}  \label{der un}
	\| \partial_\om u_n \|_{0} \leq 
C \,	\frac{\m}{\g^2 \om}  
\end{equation}
for some $C$.
\end{lemma}

\noindent
\pf
$u_1 = h_1 \in X^{(1)}$ is defined for $(\om,\m) \in A_1$ and it solves $\om^2 (h_1)_{tt} = \D h_1 + \m P_0 g$. 
Recalling Remark \ref{rem:open} and Lemma \ref{lemma:inv}, by the classical implicit function theorem it follows that $h_1$ is differentiable w.r.t. $\om$.
Differentiating  w.r.t. $\om$ gives
\begin{equation*}  
2 \om (h_1)_{tt} + L_\om [ \partial_\om h_1] = 0 .
\end{equation*}
We apply $L_\om^{-1}$ and by \eqref{stima inverso} and \eqref{h1 tt}
\[
\| \partial_\om h_1 \|_0 \leq C \, \frac{\m}{\g^2 \om}
\]
for some $C$. 

Assume that for $n \geq 1$
\begin{equation}  \label{induz der om}
\| \partial_\om h_k \|_0 \leq \bar C \, \frac{\m}{\g^2 \om} \, \exp[(-b+\t+1) \chi^k]
\quad 
\forall k = 1, \ldots , n.
\end{equation}
$h_{n+1}$ solves \eqref{comodissimo}, then it is differentiable w.r.t. $\om$. Differentiating \eqref{comodissimo} gives 
\[
F'(u_n) [ \partial_\om h_{n+1} ] = 
- 2 \om (h_{n+1})_{tt} 
+ \m f''(u_n) [\partial_\om u_n , h_{n+1} ] 
+ \m \partial_\om ( Q(u_{n-1},h_n) ) .
\]
We apply $F'(u_n)^{-1}$ and observe that
\[
\| F'(u_n)^{-1} [ \om (h_{n+1})_{tt} ] \|_0 \leq
\frac{K_1 \om }{\g} \, \|  (h_{n+1})_{tt} \|_{\t-1} 
\leq \frac{C \m}{\g^2 \om} \, \exp[(-b+\t+1)\chi^{n+1}]
\]
by \eqref{hk tt}. To estimate 
\begin{equation}  \label{termone}	
K_1 \frac{\m}{\g} \, \| f''(u_n) [\partial_\om u_n , h_{n+1} ] 
+ \partial_\om ( Q(u_{n-1},h_n) ) \|_{\t-1}
\end{equation}
we write all the integral terms and apply \eqref{S1} and \eqref{induz der om} to each of them. We write here the calculations for two terms, the other ones are analogous.
First,
\begin{equation} \label{first}
	\| \D (\partial_\om u_n) \int_\Om \gr u_n \circ \gr h_{n+1} \|_{\t-1}
\leq 
\| \partial_\om u_n \|_{\t+1} \| u_n \|_2 \| h_{n+1} \|_0 
\end{equation}
and $\| \partial_\om u_n \|_{\t+1} \leq \| \partial_\om u_n \|_0 \, N_n^{\t+1}$,
\begin{eqnarray*}
&& \| \partial_\om u_n \|_{0} 
\leq  \sum_{k= 1}^n \| \partial_\om h_k \|_{0}
\leq C \, \frac{\m}{\g^2 \om} \, 
\sum_{k\geq 1} \exp[(-b+\t+1) \chi^k] = \frac{C' \m}{\g^2 \om} \,,
\vspace{5pt}\\
&& \| u_n \|_2  \leq  \| u_n \|_{\t+1} < C \m / \g \quad \text{by \eqref{hk t+1},}
\vspace{5pt}\\
&&  N_n^{\t+1} \| h_{n+1} \|_0 
\leq 
\frac{K \m}{\g} \, \exp[(-b+\t+1) \chi^{n+1}] \quad \text{by \eqref{hn+1},}
\end{eqnarray*}
so that \eqref{first} $\leq C (\m/\g)^2 (\m/\g^2 \om) \exp[(-b+\t+1) \chi^{n+1}]$.
As second example,
\[
\| \D h_{n+1} \int_\Om \gr u_n \circ \gr (\partial_\om u_n) \|_{\t-1}
\leq 
\| h_{n+1} \|_0 \, N_{n+1}^{\t+1} \| u_n \|_2 \| \partial_\om u_n \|_0 \,.
\]
At the end we have
\[
\eqref{termone} \leq C \, \Big( \frac{\m}{\g} \Big)^3  \frac{\m}{\g^2 \om} \, \exp[(-b+\t+1) \chi^{n+1}]
\]
for some $C$, so that \eqref{induz der om} holds true for $k=n+1$ provided 
$\m/\g$ is small enough, independently on $n$. 

Finally, $\| \partial_\om u_n \|_0 \leq \sum_{k= 1}^n \| \partial_\om h_k \|_0 $ and \eqref{induz der om} implies \eqref{der un}.
\str \bs

\begin{lemma} \label{lemma:Cantor}
\textnormal{(The Cantor set).} There exist $\d_3 \leq \d_2 $ such that the Cantor set $\mathcal A_\g := A_\infty \cap \,\{ (\om,\m) : \m < \d_3 \g \}$, $\g \in (0,\lm_1)$,  has the following measure property. 

For every interval $I = \bomm$ with $0<\bom_1<\bom_2<\infty$ 
there is a constant $\bar C$ depending on $I$ such that, denoted by $\mathcal R_\g$ the rectangular region  $\mathcal R_\g = I \times (0,\d_3\g)$, 
\[
\frac{ | \mathcal R_\g \cap \mathcal A_\g |}{| \mathcal R_\g |} \, > 1 - \bar C \g \,.
\]
\end{lemma}

\noindent
\pf
We fix $\m$, we recall that $A_n(\m) := \{ \om : (\om,\m) \in A_n \}$ and define
\[
E_{n} :=  A_n(\m) \setminus A_{n+1}(\m) , \quad n \in \N .
\]
We have to prove that $\cup_{n \in \N} E_n$ has small measure. As a consequence, its complementary set $A_\infty(\m) := \cap_{n \in \N} A_n(\m)$ will be a large set. Let 
\[
\Om_{j,l}^{n} := \Big\{ \om : |  \om p_l^{(n)} (\om,\m) - \lm_j | \leq \frac{\g}{\lm_j^\t} \Big\} \,.
\]
We note that $\Om_{j,0}^{n}=\emptyset$ for all $j,n$ because $\g < \lm_1$ and $p_l^{(n)}=0$ for $l=0$.
Suppose that $\om \in \Om_{j,l}^n $. Then $C \lm_j < \om l < C' \lm_j$ for some $C,C'$ by \eqref{comparison lmm}. 
Moreover 
\[
 | \partial_\om p_l^{(n)} | \leq 2 l \m \| \partial_\om a_n \|_{\infty}
\]
by \eqref{comparison alpha beta}. 
Recalling the definition \eqref{an} of $a_n$, by \eqref{der un} and \eqref{hk t+1}
\[
\| \partial_\om a_n \|_{H^1} = \| 2 \int_\Om \gr u_n \circ \gr (\partial_\om u_n) \, dx \|_{H^1} \leq 2 \| u_n \|_2 \|  \| \partial_\om u_n \|_0  
\leq C \, \frac{\m^2}{\g^3 \om} 
\]
for some $C$, so that
\[
 | \partial_\om p_l^{(n)} | \leq C \,l \, \frac{\m^3}{\g^3 \om} \,.
\]
By \eqref{comparison lmm} it follows that
\begin{equation}  \label{vina}
 \partial_\om ( \om \, p_l^{(n)}(\om,\m) ) 
 \geq p_l^{(n)} - \om C \,l \, \frac{\m^3}{\g^3 \om}  >\, \frac{l}{4} 
\end{equation}
provided $\m/\g $ is small enough, say $\m/\g < \d_3$.
Fix  $0 < \bom_1 < \bom_2 < \infty$. 
If $\Omega_{j,l}^n \cap \bomm$ is nonempty, then 
\begin{equation}  \label{8g}
 |\Omega_{j,l}^n| < \frac{8 \g }{l \lm_j^\t} < C \bom_2 \, \frac{\g }{\lm_j^{\t+1}}\,,
 \quad
l \in \Big( \frac{C' }{\bom_2}\,\lm_j \, , \, \frac{C'' }{\bom_1} \, \lm_j  \Big) =: \Lambda(j)
\end{equation}
for some $C,C',C''$. Since $E_0 = \cup_{\lm_j \leq N_1, \, l\geq 1} \Om_{j,l}^0$,
we have the estimate
\begin{equation}  \label{E0}
|E_0 \cap \bomm | \leq \sum_{\lm_j \leq N_1} \sum_{l \in \Lambda(j)} |\Omega_{j,l}^0|
<  \g  \tilde C  \sum_{\lm_j \leq N_1} \frac{1}{\lm_j^\t} 
\end{equation}
for some $\tilde C $ depending on $\bomm$.

To estimate $|E_{n}\cap \bomm |$, $n \geq 1$, we notice that
\[
E_{n} = \bigcup_{\lm_j \leq N_{n+1},\, l \geq 1} \Omega_{j,l}^{n} \cap A_n(\m) \,.
\]
For the sets $\Om_{j,l}^{n}$ with $N_n < \lm_j \leq N_{n+1} $ we use \eqref{8g} 
and we get
\[
\Big| \bigcup_{\genfrac{}{}{0pt}{2}{N_n < \lm_j \leq N_{n+1}}{l \in \Lambda(j)}} 
 \Omega_{j,l}^{n} \cap \bomm \Big|
< 
\g \tilde C \sum_{N_n < \lm_j \leq N_{n+1}} \frac{1}{\lm_j^\t} 
\]
where $\tilde C $ is the constant of \eqref{E0}. 
To estimate the remaining sets, suppose that $\om \in \Om_{j,l}^{n} $ for some $\lm_j \leq N_n$, $l \geq 1$. 
Then by \eqref{comparison alpha beta}
\begin{eqnarray*}
| \lm_j - \om p_l^{(n-1)} | 
& \leq &
| \lm_j - \om p_l^{(n)} | + \, \om \,|p_l^{(n)} - p_l^{(n-1)} | \\
& \leq & \frac{\g}{\lm_j^\t} + 
2 \om l \m \| a_n - a_{n-1} \|_\infty \,. 
\end{eqnarray*}
Since $\om l \leq C \lm_j$ and 
$\| a_{n} - a_{n-1} \|_{H^1} \leq \| h_{n} \|_0 \| 2 u_{n-1} + h_{n} \|_2$, 
by \eqref{hk},\eqref{hk t+1} we have   
\[
| \lm_j - \om p_l^{(n-1)} | \leq \frac{\g}{\lm_j^\t} \, +
C \lm_j \frac{\m^3}{\g^2} \, \exp(-b \chi^n) .
\]
Thus 
\[
\Om_{j,l}^{n} \cap A_n(\m) \subseteq 
\Big\{ \om \,:\, 
\frac{\g}{\lm_j^\t} \, < | \lm_j - \om p_l^{(n-1)} |
\leq \frac{\g}{\lm_j^\t} \, + C \lm_j \frac{\m^3}{\g^2} \, \exp(-b \chi^n)  \Big\}
\]
and by \eqref{vina}
\[
| \Om_{j,l}^{n} \cap A_n(\m) | \leq C \bom_2 \frac{\m^3}{\g^2} \, \exp(-b \chi^n) \,.
\]
It follows that
\[
\Big| \bigcup_{\genfrac{}{}{0pt}{2}{ \lm_j \leq N_n}{l \in \Lambda(j)}} 
 \Omega_{j,l}^{n} \cap A_n(\m) \cap \bomm \Big|
\leq 
\g  C \frac{\m^3}{\g^3}  \sum_{\lm_j \leq N_n} \lm_j \exp(-b \chi^{n}) 
\]
for some $ C$ depending on $\bomm$. Now, by \eqref{Weyl} 
$\lm_j \leq N_n $ implies $ j \leq (N_n / C)^d $, then
\[
\sum_{\lm_j \leq N_n} \lm_j \,
\leq \sum_{ j \leq  (N_n/C)^d} C' j^{1/d} 
\leq C'' \int_0^{(N_n/C)^d} \xi^{1/d} d\xi 
\leq C''' N_n^{d+1} 
\]
for some $C'''$. As a consequence 
\[
\sum_{\lm_j \leq N_n} \lm_j \exp(-b \chi^{n}) 
\leq
C \exp[(-b + d +1) \chi^n] .
\]
Then
\[
|E_{n}\cap \bomm |
\leq
\g C \Big( \sum_{N_n < \lm_j \leq N_{n+1}} \frac{1}{\lm_j^\t} + 
 \exp[(-b + d +1) \chi^n] \Big) 
\]
for some $ C$, and
\[
\big| \bigcup_{n \in \N} E_n \cap \bomm \big| 
\leq 
\g C \Big( \sum_{j \geq 1} \frac{1}{\lm_j^\t} + \sum_{n \in \N} \exp[(-b + d +1) \chi^n] \Big) \,.
\]
The first series converges because by \eqref{Weyl}
\[
\sum_{j \geq 1} \frac{1}{\lm_j^\t} \leq C \sum_{j\geq 1} \frac{1}{j^{\t/d}} < \infty
\]
being $\t > d$. The second series converges because $ b  > d+1$. 
Thus 
\[
\big| \bigcup_{n \in \N} E_n \cap \bomm \big| \leq C \g 
\]
and the relative measure of $A_\infty(\mu)$ in $\bomm$ satisfies 
\[
\frac{ | A_\infty(\m) \cap \bomm | }{ \bom_2 - \bom_1 } 
\geq 1 -  C \g
\]
for some $ C$. 

Finally, we integrate in $\m$ in the interval where $\m/\g < \d_3$,
\[
| A_\infty \cap R_\g | = \int_0^{\d_3 \g} | A_\infty(\m) \cap \bomm | \, d\m \,.
\]
\str \bs

\section{\large{Inversion of the linearised operator}} \label{sec:inv}
In this section we prove Lemma \ref{lemma:inv}.
Let $u \in X^{(n)}, h \in X^{(n+1)}$. The linearised operator is 
\begin{eqnarray*}
	F'(u) h = L_\om h - \m f'(u)[h] = Dh + Sh
\end{eqnarray*}
where we split $F'(u)$ in a diagonal part 
\begin{equation*} 
	Dh := \om^2 h_{tt} - \D h \Big( 1 + \m \int_\Om | \gr u|^2 \, dx \Big) 
\end{equation*}
and a ``projection'' part
\begin{equation*} 
	Sh := - \m \D u  \int_\Om 2 \gr u \circ \gr h \, dx \, .
\end{equation*}
We recall here some results on Hill's problems. 
The proof is in the Appendix.

\begin{lemma}  
\textnormal{(Hill's problems). }
\label{lemma:Hill}
Let $\alpha (t)$ be $2\p$-periodic and $\| \a \|_\infty < 1/2$.
The eigenvalues $p^2$ of the periodic problem 
\begin{equation}  \label{rho Hill}
	\begin{cases}
	y'' + p^2 (1 + \a(t)) \, y = 0  \\
	y(t) = y(t + 2\p)
	\end{cases}
\end{equation}
form a sequence $\{p_l^2  \}_{l \in \N }$ such that 
\begin{equation}  \label{comparison lmm}
	\frac13 \, l \leq p_l \leq 2l
	\quad \, \forall l \in \N .
\end{equation}
For $\a,\b $ $2\p$-periodic, $\| \a \|_\infty, \| \b \|_\infty < 1/2$,
\begin{equation}  \label{comparison alpha beta}
\big| p_l(\alpha) - p_l(\beta) \big| \leq 2 l \, \| \alpha - \beta \|_{\infty} 
\quad \, \forall l \in \N .
\end{equation}
 
The eigenfunctions $\psi_l(t)$ of \eqref{rho Hill} form an orthonormal basis 
of $L^2(\T)$ w.r.t. the scalar product
\begin{equation*}  
(u,v)_{L^2_\a} = \int_0^{2\p} u v (1+\a) \, dt 
\end{equation*}
and also an orthogonal basis of $H^1(\T)$ w.r.t. the scalar product
\[
(u,v)_{H^1_\a} = \int_0^{2\p} u'v' \, dt + (u,v)_{L^2_\a}.
\]
The corresponding norms are equivalent to the standard Sobolev norms,
\begin{equation}  \label{norme equivalenti}
\frac12 \, \| y \|_{L^2} \leq \|y \|_{L^2_\a} \leq 2 \, \| y \|_{L^2},
\quad 
\frac12 \, \| y \|_{H^1} \leq \|y \|_{H^1_\a} \leq 2 \, \| y \|_{H^1}  .
\end{equation}
\end{lemma}

\begin{lemma}  \label{inv D}
\textnormal{(Inversion of $D$). }
Let $u \in X^{(n)}$, $a(t) := \int_\Om | \gr u|^2 \, dx \in H^1(\T)$, 
$ \| a \|_{H^1} < 1$, 
$ \| a \|_{\infty} < 1/2$. 
Let $p_l^2$ be the eigenvalues of the Hill's problem
\begin{equation}  \label{Hill D}
	\begin{cases}
	y'' + p^2 (1 + \m a(t) ) \, y = 0 \\
	y(t) = y(t + 2\p) \,.
	\end{cases}
\end{equation}
If $(\om,\m)$ satisfy the non-resonant condition
\begin{equation*}  
| \om p_l - \lm_j| > \frac{\g }{\lm_j^{\t}} 
\quad 
\forall \lm_j \leq N_{n+1}, \ \  l \in \N \, ,
\end{equation*} 
then $D$ is invertible, $D^{-1} : X^{(n+1)} \to X^{(n+1)}$ and
\begin{equation}  \label{D -1}
\| D^{-1} h \|_{0} 
\leq \frac{C }{\g  } \,  \| h \|_{\t - 1}
\quad \forall h \in X^{(n+1)}  
\end{equation}
for some constant $C$. 
\end{lemma}

\noindent
\pf
If $h = \sum h_j(t) \, \ph_j(x) $, then $Dh = \sum D_j h_j(t) \, \ph_j(x) $, 
where
\begin{equation*} 
D_j z(t) = \om^2 z''(t) + \lm_j^2 z(t) \rho(t) ,
\quad \rho(t) := 1 + \m a(t).
\end{equation*}
Using the eigenfunctions $\psi_l(t)$ of \eqref{Hill D} as a basis of $H^1(\T)$,
\[
D_j z(t) 
= \sum_{l \in \N} (\lm_j^2 - \om^2 p_l^2 ) \hat z_l \,\psi_l(t) \,\rho(t)  ,
\qquad
z= \sum_{l \in \N} \hat z_l \, \psi_l(t)  ,
\]
and $K_j := (1/\rho) D_j$ is the diagonal operator 
$\{ \lm_j^2 - \om^2 p_l^2 \}_{l \in \N}$.
Since  $|\lm_j^2 - \om^2 p_l^2 | >  \g / \lm_j^{\t -1}$ for all $\lm_j \leq N_{n+1}$,
we have that $K_j$ is invertible and 
\[
\| K_j^{-1} z \|_{H^1_{\m a}}^2 = 
\sum_{l \in \N} \Big( \frac{\hat z_l}{\lm_j^2 - \om^2 p_l^2} \Big)^2  
\| \psi_l \|_{H^1_{\m a}}^2 
\leq 
\frac{ \lm_j^{2(\t-1)} }{ \g ^2 } \, \| z \|^2_{H^1_{\m a}} \, .
\]
By \eqref{norme equivalenti} 
$ \| K_j^{-1} z \|_{H^1} \leq 4 (\lm_j^{\t-1} / \g ) \, \| z \|_{H^1} $. 
Since $D_j^{-1} z = K_j^{-1} (z / \rho)$ 
and $\| 1/ \rho \|_{H^1}$ is smaller than a universal constant, 
\begin{equation*}  
\| D_j^{-1} z \|_{H^1} \leq  \frac{C \lm_j^{\t-1} }{\g } \, \| z \|_{H^1} \,.
\end{equation*}
Since $D^{-1} h = \sum_j D^{-1}_j h_j(t) \ph_j(x)$ we obtain \eqref{D -1}. 
\str \bs

\smallskip

\begin{lemma} \label{lemma:S} 
\textnormal{(Control of $S$). }
For all $s \geq 0$, if $u \in X_{s+2}$ then $S : X_{0} \to X_{s}$ is bounded and
\begin{equation*}
\| Sh \|_{s} \leq \mu \| u \|_{s+2} \| u \|_{2} \| h \|_{0}
\qquad
\forall h \in X_0  \,.
\end{equation*}
\end{lemma}

\noindent
\pf
Since $\int_\Om \gr u \circ \gr h \, dx$ does not depend on $x$, 
\[
\Big\| \D u  \int_\Om \gr u \circ \gr h \, dx \Big\|_{s}
\leq \| \D u \|_{s} \, \Big\| \int_\Om \gr u \circ \gr h \, dx \Big\|_{H^1(\T)} \, .
\]
$\int_\Om \gr u \circ \gr h \, dx = \sum_j \lm_j^2 u_j(t) h_j(t)$, 
so 
$\| \int_\Om \gr u \circ \gr h \, dx\|_{H^1(\T)}$ $\leq$ $ \| u \|_{2} \| h \|_{0}$ 
by H$\ddot{\mathrm{o}}$lder inequality. 
\str \bs

\bigskip

\noindent
\textbf{Proof of Lemma \ref{lemma:inv}. }
$F'(u) = D + S = ( I + S D^{-1}) D$ 
where $I$ is the identity map. Since $D^{-1}$ satisfies \eqref{D -1}, 
we have to prove the invertibility of $I + S D^{-1}$ in norm $\| \ \|_{\t-1}$.
By Neumann series it is sufficient to show that
\begin{equation}  \label{SD-1}
	\| S D^{-1} h \|_{\t-1} \leq \frac12 \, \| h \|_{\t-1} 
	\quad
	\forall h \in X^{(n+1)} .
\end{equation}
By Lemmas \ref{inv D} and \ref{lemma:S} 
\[
\| S D^{-1} h \|_{\t-1} 
\leq \m \| u \|_{\t-1+2} \| u \|_{2} \| D^{-1} h \|_{0} 
\leq
 \frac{C\m}{\g }\, \| u \|_{\t+1}^2 \| h \|_{\t-1} 
\]
because  $\| u \|_2 \leq  \| u \|_{\t+1} $. 
Thus the condition 
\[
\frac{\m}{\g}\, \| u \|_{\t+1}^2 \leq \frac1{2C} \, =: K_1'
\]
implies \eqref{SD-1} and by Neumann series 
$\| ( I + S D^{-1})^{-1} h \|_{\t-1}  \leq 2 \| h \|_{\t-1}$.
\str \bs 

\smallskip 

\section{\large{Proof of the theorems}}  \label{sec:proof}

\textbf{Proof of Theorem \ref{thm:uni}.}
Let $g \in X_{\s,s_0}$ and $2d < 2(s_1 -1) < s_0$ as assumed in the theorem.
We apply Lemma \ref{lemma:indutt} with
\begin{equation*}  \label{ss1}
\t := s_1 -1 \,.	
\end{equation*}
The construction of the sequence $(u_n)$ is possible provided the parameters $(\om,\m)$ belong to $A_n$ for all $n \in \N$. 
Lemma \ref{lemma:Cantor} assures that, for $\m/\g$ sufficiently small, the set $\mA_\g$ of parameters satisfying this property is a nonempty set, which is very large in a Lebesgue measure sense.
Lemmas \ref{lemma:existence} and \ref{lemma:uniqueness} complete the proof.

\str \bs

\jump

\noindent
\textbf{Proof of Theorem \ref{thm:periodic}.} Let $g \in X_{\s,s_0}$ with $2d < 2(s_1 -1) < s_0$ as assumed in the theorem.
We consider a Lyapunov-Schmidt reduction splitting the space $\tilde X_{\s,s}$ in two subspaces $\tilde X_{\s,s} = Y \oplus (W \cap \tilde X_{\s,s})$,
\[
Y := \{ y(t) \in H^1(\T,\R) \} \,,
\quad
W := \Big\{ w \in \tilde X_{0,0} \, : \, 
w(x,t) = \sum_{j \geq 1} w_j(t) \tilde \ph_j(x) \Big\} \,.
\]
We denote $\Pi_Y, \Pi_W$ the projectors on $Y,W$, and observe that $\Pi_Y$ is the map
\[
u \mapsto \int_{(0,2\p)^d} u(x,t) \, dx .
\]
We define 
\[
g_0 (t) := \Pi_Y g, 
\quad
\bar g (x,t) := \Pi_W g .
\]
We decompose $u(x,t) = y(t) + w(x,t)$, $y \in Y$, $w \in W$, and note that 
\[
f(u) = f(y+w) = f(w) \in W.
\]
Then projecting equation \eqref{K} on $Y$ gives 
\begin{equation} \label{Y-eq}
\om^2 y''(t) = \m g_0 (t)	\quad (Y \ equation)
\end{equation}
while projecting it on $W$ 
\begin{equation} \label{W-eq}
 L_\om w = \m ( f(w) + \bar g) \quad (W \ equation).
\end{equation}
Equation \eqref{Y-eq} is an ODE. 
With direct calculations (or by Fourier series) we see that \eqref{Y-eq} admits $2\p$-periodic solutions if and only if 
\begin{equation} \label{just}
\int_0^{2\p} g_0(t) \, dt = 0 
\end{equation}
and \eqref{just} is just assumption \eqref{hyp:ort}. We note that, if $y(t)$ solves 
\eqref{Y-eq}, then also $y(t) + c$ solves \eqref{Y-eq}, for all $c \in \R$.
Moreover, the unique solution $y(t)$ of \eqref{Y-eq} such that $\int_0^{2\p} y(t) \, dt = 0$ satisfies
\[
\| y \|_{H^1}  \leq 
\| y'' \|_{H^1} \leq \frac\m {\om^2} \, \| g_0 \|_{H^1}.
\] 

To solve \eqref{W-eq}, we consider all the calculations in Sections 
\ref{sec:iteration},\ref{sec:solution},\ref{sec:Cantor},\ref{sec:inv} replacing
 $X_{\s,s}$ with $\tilde X_{\s,s} \cap W$ and 
 $\lm_j, \ph_j(x)$ with $\tilde \lm_j, \tilde \ph_j(x)$, $j \geq 1$. 
It follows the existence of a unique solution $w \in \tilde  X_{\s,s_1} \cap W $ of \eqref{W-eq} satisfying 
\[
\| w \|_{\s,s_1} \leq \frac{\m}{\g } \, C ,
\quad
\| w_{tt} \|_{\s,s_1-2} \leq \frac{\m}{\g \om^2} \, C .
\]
Then $u = y + w$ solves \eqref{K}\eqref{perbc}. Since 
\[
\| u \|_{\s,s_1}^2 = \| y \|_{H^1}^2 + \| w \|_{\s,s_1}^2 ,
\quad
\| u_{tt} \|_{\s,s_1-2}^2 = \| y'' \|_{H^1}^2 + \| w_{tt} \|_{\s,s_1-2}^2 ,
\]
we obtain estimates \eqref{estim sol per}.
\str \bs

\section{\large{Appendix}}

\textbf{Proof of Lemma \ref{lemma:Hill}.}
The proof follows from classical results in \cite{E,Kato}. 
First, if $y'' + p^2 (1+\a) y = 0$, then
\[
\int_{0}^{2\p} y'^2 \, dt = p^2 \int_0^{2\p} (1+\a) y^2 \, dt \, ,
\]
so that $p^2 \geq 0$ because $ 1+\a $ is positive. $p_0^2=0$ is an eigenvalue, the corresponding 
eigenfunctions are the constants, and all the other eigenvalues are positive.

By \cite[Theorem 2.2.2, p.23]{E}, for every $k \in \N$ 
both $p_{2k+1}$ and $p_{2k+2}$ satisfy
\[
\frac23 (k+1)^2 \leq p^2 \leq 2 (k+1)^2 
\]
and \eqref{comparison lmm} follows. 

The equivalence \eqref{norme equivalenti} and 
the orthogonality of $(\psi_l)$ w.r.t.\:$( \, , \, )_{H^1_\a}$ 
can be verified by direct calculations.

To prove \eqref{comparison alpha beta}, we define 
\[
q(\th)(t) = q (t) := 1 + \a(t) + \th (\beta(t) - \a(t)) ,
\quad \th \in [0,1]
\]
and the ``Liouville's change of variable'' $t \to \xi$
\[
t = f(\xi) \quad \Leftrightarrow \quad \xi = g(t) := \dfrac{1}{c} \, \int_0^{t} \sqrt{q(s)} \, ds,
\qquad
c :=  \dfrac{1}{2\p} \, \int_0^{2\p} \sqrt{q(t)} \, dt.
\]
We note that $p^2, y(t)$ satisfy
\begin{equation} \label{Hill yt}
\begin{cases}
y''(t) + p^2 q(t) y(t) =0 \\
y(t ) = y(t+ 2 \p)
\end{cases}
\end{equation}
if and only if $p^2, z(\xi)$ satisfy
\begin{equation} \label{Hill zxi}
\begin{cases}
z''(\xi) +  c^2 [ p^2 - Q(f(\xi))] z(\xi)  =0 \\
z(\xi ) = z(\xi+ 2 \p)
\end{cases}
\end{equation}
where 
\[
z(\xi) := y(f(\xi)) \sqrt[4]{ q(f(\xi)) } ,
\quad
Q(t) := - \frac{5}{16} \frac{q'(t)^2}{q(t)^3} \, + \frac14 \frac{q''(t)}{q(t)^2} \,.
\]
The  operators $T(\th) : z \mapsto -z'' + c^2 Q(f) z$ are selfadjoint in $L^2(0,2\p)$.  
We apply \cite[Theorem 3.9, VII.3.5, p.\,392]{Kato} 
to the holomorphic family $\{ T(\th) : \th \in [0,1] \}$ (see \cite[Definition VII.2.1, p.\,375, Example 2.12, VII.2.3, p.\,380 and Example 6.13, III.6.8, p.\,187]{Kato} to verify the hypotheses of the Theorem in the present case) 
to prove that the eigenvalues and eigenfunctions of \eqref{Hill zxi} are analytic in $\th$.  As a consequence, the eigenvalues and eigenfunctions of \eqref{Hill yt} are analytic in $\th$ as well.
This allows us to differentiate the equation 
\[
\psi_l(\th)'' + p_l(\th)^2 (1 + \a + \th(\b - \a)) \psi_l(\th) = 0
\]
w.r.t.\,$\th$. 
Recalling that $\int_0^{2\p} q \psi_l^2 \, dt = 1$, multiplying by $\psi_l(\th)$ and integrating 
\begin{equation*}  
\partial_\th p_l(\th) = - \frac12\, p_l(\th) 
\int_0^{2\p} 
(\b - \a) \psi_l(\th)^2 \, dt \, .
\end{equation*}
Since $p_l(\th) \leq 2l$ and $q(\th) \geq 1/2$,
\[
|\partial_\th p_l(\th)| \leq l \int_0^{2\p} \frac{ |\b - \a| }{ q(\th)} \, q(\th)
\psi_l(\th)^2 \, dt \leq 2 l \| \b - \a \|_\infty
\]
and
\[
|p_l(\b) - p_l(\a) | \leq  
\int_0^1 | \partial_\th p_l(\th) | \, d\th \leq 2 l \| \b - \a \|_\infty.
\]
\str \bs

\jump

\noindent
\textbf{Acknowledgements.} The author would like to thank Massimiliano Berti for his encouragement and many suggestions and John Toland for very useful discussions.

\small{

\end{document}
\begin{thebibliography}{99}   

\bibitem{Ar-81} A. Arosio, \textit{Asymptotic behaviour as $t \to +\infty$ of the solutions of linear hyperbolic equations with coefficients discontinuous in time (on a bounded domain)}, J. Diff. Eq. \textbf{39} (1981), n.2, 291--309.

\bibjump

\bibitem{Ar} A. Arosio, \textit{Averaged evolution equations. The Kirchhoff string and its treatment in scales of Banach spaces}, in: 2nd Workshop on functional-analytic methods in complex analysis (Trieste, 1993), World Scientific, Singapore.

\bibjump

\bibitem{AS-84} A. Arosio, S. Spagnolo, \textit{Global solutions of the Cauchy problem for a nonlinear hyperbolic equation}, in: Nonlinear PDE's and their applications,   Coll\'ege de France Seminar, Vol. VI, 1--26, H. Brezis \& J.L. Lions eds., Research Notes Math. \textbf{109}, Pitman, Boston, 1984.

\bibjump

\bibitem{AP-96} A. Arosio, S. Panizzi, \textit{On the well-posedness of the Kirchhoff string}, Trans. Amer. Math. Soc. \textbf{348} (1996), no.1, 305--330.
  
\bibjump

\bibitem{BaBe-Forced} P. Baldi, M. Berti, \textit{Forced vibrations of a nonhomogeneous string}, preprint, 2006.  
 
\bibjump

\bibitem{Bernstein} S.N. Bernstein, \textit{Sur une classe d'\'equations fonctionelles aux d\'eriv\'ees partielles}, Izv. Akad. Nauk SSSR Ser. Mat. \textbf{4} (1940), 17--26. 
  
\bibjump

\bibitem{Be-Topics} M. Berti, \textit{Topics on ``Nonlinear oscillations of Hamiltonian PDEs''}, to appear in
Progress in nonlinear differential equations and their applications, Birkh\"auser,
Boston.

\bibjump

\bibitem{BeBo-Ck} M. Berti, P. Bolle, \textit{Cantor families of periodic solutions of wave equations with $C^k$ 
nonlinearities}, preprint 2007.

\bibjump

\bibitem{BeBo-CantorFam} M. Berti, P. Bolle, \textit{Cantor families of periodic solutions for completely resonant nonlinear wave equations}, Duke Math. J. \textbf{134} (2006), no.\,2, 359--419.
 
\bibjump

\bibitem{Bou-IMRN-1994} J. Bourgain, \textit{Construction of quasi-periodic solutions for Hamiltonian perturbations of linear equations and applications to nonlinear PDE}, Int. Math. Res. Notices \textbf{11} (1994), 475--497.

\bibjump

\bibitem{Bou-GAFA-95} J. Bourgain, \textit{Construction of periodic solutions of nonlinear wave equations in higher dimension}, Geom. and Funct. Anal. \textbf{5} (1995), 629--639.

\bibjump

\bibitem{Bou-Annales-98} J. Bourgain, \textit{Quasi-periodic solutions of Hamiltonian perturbation of 2D linear Schr{\accent 127o}dinger equations}, Ann. of Math. \textbf{148} (1998), 363--439.

\bibjump

\bibitem{Bou-Chicago} J. Bourgain, \textit{Periodic solutions of nonlinear wave equations}, Harmonic analysis and partial differential equations, 69--97, Chicago Lectures in Math., Univ. Chicago Press, Chicago, IL, 1999.

\bibjump

\bibitem{Car} G.F. Carrier, \textit{On the nonlinear vibration problem of the elastic string}, Quart. Appl. Math. \textbf{3} (1945), 157--165; 
--------, \textit{A note on the vibrating string}, Quart. Appl. Math. \textbf{7} (1949), 97--101.

\bibjump

\bibitem{Petits} W. Craig,  {\it Probl\`emes de petits diviseurs dans 
les \'equations aux d\'eriv\'ees partielles}, 
Panoramas et Synth\`eses, 9, Soci\'et\'e Math\'ematique de France, Paris, 2000.

\bibjump

\bibitem{CW} W. Craig, E. Wayne, 
\textit{Newton's method and periodic solutions of nonlinear wave equations}, 
Comm. Pure Appl. Math. \textbf{46} (1993), 1409--1501.

\bibjump

\bibitem{DS-InvMath-92} P. D'Ancona, S. Spagnolo, \textit{Global solvability for the degenerate Kirchhoff equation with real analytic data}, Invent. Math. \textbf{108} (1992), 247--262.

\bibjump

\bibitem{Di-69} R.W. Dickey, \textit{Infinite systems of nonlinear oscillation equations related to the string}, Proc. Amer. Math. Soc. \textbf{23} (1969), no.3, 459--468.

\bibjump

\bibitem{E} M.S.P. Eastham, \textit{The spectral theory of periodic differential equations}, Scottish Academic Press Ltd., Edinburgh, 1973.
  
\bibjump

\bibitem{Fu} D. Fujiwara, \textit{Concrete characterization of the domains of fractional powers of some elliptic differential operators of the second order}, Proc. Japan Acad. \textbf{43} (1967), 82--86.

\bibjump

\bibitem{Lax} L. Glimm, P. Lax, \textit{Decay of solutions of systems of hyperbolic conservation laws}, Mem. Amer. Math. Soc. \textbf{101}, AMS, Providence RI, 1970.

\bibjump

\bibitem{Kato} T. Kato, \textit{Perturbation theory for linear operators}, 2nd ed., Grundlehren der mathematischen Wissenschaften 132, Springer-Verlag, Berlin, 1976. 

\bibjump

\bibitem{K} G. Kirchhoff, \textit{Vorlesungen $\ddot{u}$ber mathematische Physik: Mechanik}, ch.29, Teubner, Leipzig, 1876.

\bibjump

\bibitem{Kl-Maj} S. Klainerman, A. Majda, \textit{Formation of singularities for wave equations including the nonlinear vibrating string}, Comm. Pure Appl. Math. \textbf{33} (1980), 241--263.

\bibjump

\bibitem{Kuk-93} S. Kuksin, \textit{Nearly integrable infinite-dimensional Hamiltonian systems}, Lecture Notes in Math. \textbf{1556}, Springer, Berlin, 1993.

\bibjump

\bibitem{JLL-78} J.L. Lions, \textit{On some questions in boundary value problems of mathematical physics}, in: Contemporary developments in continuum mechanics and PDE's, G.M. de la Penha \& L.A. Medeiros eds., North-Holland, Amsterdam, 1978.

\bibjump

\bibitem{LM} J.L. Lions, E. Magenes, \textit{Espaces de fonctions et distributions du type de Gevrey et problemes aux limites paraboliques}, Ann. Mat. Pura Appl. \textbf{68} (1965), 341--417; 
--------, \textit{Problemes aux limites non homogenes et applications}, Dunod, Paris, 1968.

\bibjump

\bibitem{Man-05} R. Manfrin, \textit{On the global solvability of 
Kirchhoff equation for non-analytic initial data}, J. Differential Equations \textbf{211} (2005), 38--60.

\bibjump

\bibitem{Moser-61} J. Moser, \textit{A new technique for the construction of solutions of nonlinear differential equations}, Proc. Nat. Acad. Sci. U.S.A. \textbf{47} (1961), 1824--1831.

\bibjump

\bibitem{Moser-Pisa66} J. Moser, \textit{A rapidly convergent iteration method and non-linear differential equations. I. II.}, Ann. Scuola Norm. Sup. Pisa (3) \textbf{20} (1966), 265--315 \& 499--535. 

\bibjump

\bibitem{Na} R. Narasimha, \textit{Nonlinear vibration of an elastic string}, J. Sound Vibration \textbf{8} (1968), 134--146.

\bibjump

\bibitem{Pl-Tol} P. Plotnikov, J. Toland, \textit{Nash-Moser theory for standing water waves}, Arch. Ration. Mech. Anal. \textbf{159} (2001), 1--83.

\bibjump

\bibitem{Pok-75} S.I. Pokhozhaev, \textit{On a class of quasilinear hyperbolic equations}, Mat. Sbornik \textbf{96} (1975), 152--166 (English transl.: Mat. USSR Sbornik \textbf{25} (1975), 145--158).

\bibjump

\bibitem{Poschel-Pisa} J. P\accent 127oschel, \textit{A KAM theorem for some nonlinear PDEs}, Ann. Scuola Norm. Sup. Pisa, cl. Sci., IV ser. \textbf{15} (1996), n.\,23, 119--148.

\bibjump

\bibitem{Rab} P. Rabinowitz, \textit{Periodic solutions of nonlinear hyperbolic partial differential equations. II}, Comm. Pure Appl. Math. \textbf{22} (1969), 15--39.

\bibjump

\bibitem{RS} M. Reed, B. Simon, \textit{Methods of modern mathematical physics}, Academic Press, Inc., New York, 1978.

\bibjump

\bibitem{Sp-Milano-94} S. Spagnolo, \textit{The Cauchy problem for Kirchhoff equations}, Rend. Sem. Mat. Fis. Milano \textbf{62} (1994), 17--51.

\bibjump

\bibitem{W-90} C.E. Wayne, \textit{Periodic and quasi-periodic solutions of nonlinear wave equations via KAM theory}, Comm. Math. Phys. \textbf{127} (1990), 479--528.

\bibjump

\bibitem{Z} E. Zehnder, \textit{Generalized Implicit Function Theorems}, Chapter VI, in: L. Nirenberg, \textit{Topics in nonlinear functional analysis}, Courant Inst. of Math. Sciences, New York Univ., New York, 1974.


\end{thebibliography}
